\documentclass[11pt]{article}
\usepackage{amsmath,amssymb}
\usepackage{enumerate}
\usepackage{color}

\usepackage{xy}
\usepackage{xypic}
\usepackage{subfigure}
\usepackage{tikz-cd}

\newcommand{\R}{\mathbb{R}}

\newcommand{\N}{\mathbb{N}}
\newcommand{\bea}{\begin{eqnarray}}
\newcommand{\eea}{\end{eqnarray}}
\newcommand{\la}{\lambda}

\def\de{\delta}
\def\e{\varepsilon}

\def\s{\sigma}

\def\x{\otimes}

\def\supp{{\rm supp}}
\def\diam{{\rm diam}}

\def\1{\rm Id}
\def\sgn{\rm sgn}
\def\sup{{\rm sup}}
\def\vol{{\rm vol}}

\def\ci{\circ}
\def\cl{{\rm cl}}

\newcommand{\qed}{$\hfill\blacksquare$}

\def\V{\noindent}

\def\Int{{\rm int}}
\def\cl{{\rm cl}}

\newcommand{\bean}{\begin{eqnarray*}}
\newcommand{\eean}{\end{eqnarray*}}

\newtheorem{Theorem}{Theorem}

\newcommand{\ben}{\begin{enumerate}}
\newcommand{\een}{\end{enumerate}}
\newcommand{\bit}{\begin{itemize}}
\newcommand{\eit}{\end{itemize}}
\newcommand{\edoc}{\end{document}}

\usepackage{hyperref}

\usepackage{enumerate}

\title{On the Hauptvermutung of Causal Set Theory}

\begin{document}

\author{Olaf M\"uller\footnote{Institut f\"ur Mathematik, Humboldt-Universit\"at zu Berlin, Unter den Linden 6, D-10099 Berlin, \texttt{Email: o.mueller@hu-berlin.de}}}

\date{\today}
\maketitle

\begin{abstract}
\V We formulate the Hauptvermutung of Causal Set Theory in two mathematically well-defined but different ways one of which turns out to be wrong and the other one turns out to be true. A further result is that the Hauptvermutung is true if we replace finite with countable sets.
\end{abstract}

The idea of Causal Set Theory \cite{BLMS}, \cite{Surya}, \cite{DS}, \cite{BL}, is to replace Lorentzian geometries with order relations on subsets of $\N$. In particular, if the subset is taken to be finite, the hope is that a hypothetical integral over all Lorentzian geometries (a gravitational path integral) could be replaced with a finite sum. The Hauptvermutung (main conjecture) of Causal Set Theory (in a not entirely mathematically precise formulation) is that the spacetime should be reconstructible from the finite data if they arise from some appropriate statistical process, see e.g. \cite{Surya}, \cite{Bombelli}. The present article tries to elaborate some related well-defined statements and check their validity. The following results very likely still do not correspond exactly to what in Causal Set Theory is expected to be a solution but should be considered as an honest intent to make some version of the Hauptvermutung rigorous and mathematically accessible. Probably other versions will arise in subsequent debates within the respective scientific communities.

\bigskip

In our context, a {\bf Cauchy slab} is a Lorentzian spacetime $(X,g)$ such there is an isometric embedding  $f: (X,g) \rightarrow (N,h)$ into a globally hyperbolic spacetime $(N,h)$ with $f(X) = I^+(S_-) \cap I^- (S_+)$ for two connected spacelike Cauchy hypersurfaces $S_- $ and $S_+ \subset I^+(S_-)$ of $(N,h)$. \footnote{In previous publications, Cauchy slabs have been defined as containing the future and past boundary of the image of $f$, but here we want to exclude them for technical reasons.} A Cauchy slab is called {\bf spatially compact (s.c.)} iff the closure of $f(X)$ is compact, or equivalently, if $S_+$ (and thus also $S_-$) is compact.     

\bigskip

Let CS be the set of time-oriented isometry classes of normalized (i.e., unit-volume) spatially compact Cauchy slabs.

For $X \in CS$, define $V(X) $ as the set of those $a: \N\rightarrow X$ that satisfiy the {\em volume law} 

\bea
\label{Volume-Law}
V (X) := \{ a: \N\rightarrow X |   \forall (p,q ) \in X^2: \# (a^{-1} (J(p,q)) \cap \N_n)/n \rightarrow_{n \rightarrow \infty}  \vol_g (J(p,q)) . 
\eea

For all $a \in V(X)$, the set $a(\N)$ is everywhere dense in $X$, as open subsets have nonzero measure.

\begin{Theorem}
$\forall X \in  CS : V(X) \neq \emptyset$. 
\end{Theorem}

\V {\bf Proof.} Let $t$ be a Cauchy temporal function $t$ adapted to the boundary in the sense of $ t(S_\pm) = \{ \pm 1\}$ (existence is shown in \cite{oM-inv}). The {\bf flip metric to $g = -u dt^2 + g_t$ and $t$} is $g^t := udt^2 + g_t$; it has the same volume form as $g$. Let $D_+$ be the metric space metric derived from $g^t$. Recall the definition of Hausdorff volume: For a subset $A \subset X$ of a metric space $X$ and for $d > 0 \leq r$ we define

\[ H_{r,d} (A) := \inf \{\sum_{k=0}^\infty \diam (U_k)^d | U: \N\rightarrow 2^X, A \subset \bigcup_{k=0}^\infty  U_k, \diam (U_k) < r \forall k \in \N \}  \]

As $r \mapsto H_{r,d} (A) $ is decreasing, $H_{r,d}(A) \to_{r \to 0}  \sup \{ H_{r,d} (A) | r >0 \} =: H_d(A) \in [0; \infty]$.

$H_d$ is a metric outer measure, in particular, all Borel sets are measurable. It turns out by elementary calculations that there is one $d \in [0; \infty]$ called the {\bf Hausdorff dimension of $A$} s.t. $H_k(A) =0 \forall k>d $ and $H_k(A) = \infty \forall k<d$. Now focus on a causal diamond $J(p,q)$. We know that $J(p,q) = \cl (I(p,q))$ with $I(p,q) $ open, and that $J(p,q)$ is rectifiable: Prop. 6.3.1 of \cite{HE} implies that $J^+(p) $ and $J^-(q) $ each have a Lipschitz boundary, and each interior complementary ('null normal') vector field of $J^+(p)$ is future whereas each interior complementary vector field of $J^-(q)$ is past, so the Lipschitz submanifolds intersect each other strictly Transversally, therefore Lipschitzness holds even at their intersection (as the maximum of two Lipschitz functions is Lipschitz). So the $n$-dimensional Lebesgue volume of the boundary of $J(p,q)$ vanishes.

For $n \in \N$, let $\N_n:= \{ m \in \N | m \leq n \}$. A classical result --- the generalized 
Glivenko-Cantelli lemma, see \cite{Parthasarathy}, Th. 7.1 --- asserts that for a separable metric probability space $(X, d \mu)$ the (metrizable) product topology and the (Borel, probability) product measure $\la$ on the space $X^\N$ of sequences in $X$, there is a subset $U$ of full $\la$-measure such that for all $a \in U$ the sequence of measures
	
	\[ n \mapsto \frac{1}{n+1} \cdot (a|_{\N_n})_* (\#) ,   \]
	
assigning to a Borel subset $A$ the number $\frac{1}{n+1} \# (a^{-1} (A) \cap \N_n) = \frac{1}{n+1}\sum_{k=0}^n \delta_{a(k)}$, where $\delta_x$ is the delta distribution at $x \in X$, converges weakly to $\mu$. We conclude the proof via the Portmanteau theorem: Weak convergence of measures $\mu_n$ on a sigma-compact measure space to a limit measure $\mu_\infty$ implies for any set $A$ with $\mu_\infty( \partial A) = 0$ (as for $A= J(p,q)$) that $\mu_n (A) \to_{n \to \infty} \mu_\infty (A)$, therefore $V(X) \supset U$ is of full $\la$-measure, thus nonempty. \hfill \qed

 \newpage

For a set $X$, let $P(S)$ denote the set of subsets of $S$. Let ${\rm Or} (\N) \subset P(\N\times \N)$ be the set of order relations on $\N$. We define a map $C: CS \to P({\rm Or}(\N))$ assigning to each $X \in {\rm CS}$ the set $\{ (a,a)^{-1}(\leq) | a \in V(X) \}$, and also regard $J$ as a non-right-unique {\em relation} from $CS$ to ${\rm Or} (\N)$. This relation is never right-unique: For each order $U$ on $\N$ with $XCU$, we get $XC(\s^* U)$ for each $\s: \N \rightarrow \N$ bijective with $\# (\supp (\1 - \s))< \infty$, as $a \ci \s \in V(X)$ for all $a \in V(X)$, so we can choose $\s $ to be the elementary permutation of two elements $a(n)$ and $a(m) \geq a(n)$. 

\begin{Theorem}[The Countable Hauptvermutung is true]
\label{Countable-CausalSet}
The relation $C$ from $ CS$ to $ {\rm Or} (\N)$ is left-unique (injective), i.e., for each $U \in P(\N\times \N)$ there is at most one $X \in CS$ with $(X,U) \in C$. 
\end{Theorem}

{\bf Proof.} Each $U = a^* (\leq) = (a,a)^{-1} (\leq_X)$ is an order relation on $\N$. For $A \subset \N$ put $J_U^-(A) := \{ m \in \N | mUa \forall a \in A\}$. Let $K$ be the set of $(\leq, U)$-increasing maps\footnote{We could include partial maps, of finite and infinite domain of definition, but this makes no difference.} $\phi: \N \to \N$ s.t. $\phi(\N)$ is bounded above. We define $L:= K/\sim$ where $k \sim l :\Leftrightarrow J^+_U(k(\N)) = J^+_U(l(\N))  \land J^-_U(k(\N)) = J^-_U(l(\N))$ for all $k, l \in K$. Then two subsequences are equivalent if and only if they converge (always from below by monotonicity) to the same point of $X$. On the other hand, as $a(\N)$ is everywhere dense, for each $q \in X$ and each $\e >0$ there is $k \in \N$ such that $a(k) \in I^-(q) \cap B(q, \e) $, thus $L$ is in bijection to $X$. On $L $ we define an order relation $\leq$ by $ l_1 \leq l_2 : \Rightarrow J^-(l_1) \subset J^-(l_2)$, and a relation $\ll = \beta (\leq)$ on $L$:

$$(x,y) \in \beta (\leq) : \Leftrightarrow \big( x \leq y \land ( \exists u,v \in X:  x < u< v < y \land J^+(u) \cap J^-(v) {\rm \ not \ totally \ ordered} ) \big)$$   

We denote $I^\pm (x) := \{ y \in L | x (\beta (\leq))^{\pm 1} y \}$ and give $L$ the topology generated by the timelike diamonds $I(u,v) := I^+(u) \cap I^-(v)$ for $u,v \in L$ as a subbasis. Define a Borel measure $m$ on $L$ by 

$m (I(u,v)) :=   \lim_{n \rightarrow \infty} \frac{1}{n} \#  \{ k \in \N_n | uUa(k)Uv \}$. Therefore, the bijection between $X $ and $L$ preserves the chronological and the causal relation as well as the measure. By the Malament-Hawking-King-McCarthy Theorem \cite{M}, \cite{HKM}, any other Cauchy slab inducing the same structures on $L$ is isometric to $X$, which shows injectivity of $C$. \hfill \qed

\medskip

So we could formulate general relativity (in the s.c. case) as a theory on ${\rm Or} (\N) \subset P(\N \times \N)$.

\bigskip

As mentioned above, Causal Set Theory prefers a {\em finite} number $K$ (typically taken to be around $10^{240}$) of points. In that case, we can only hope to {\em approximately} reconstruct the spacetime, if at all. This calls for the notion of a distance between spacetimes. There is a notion of Gromov-Hausdorff distance of spacetimes we can use, defined in \cite{oM:LGH} and independently in \cite{eMsS}. It is the following: Given two Cauchy slabs $(M,g)$ and $(N,h)$ we define $d^- (M,N) := \inf \{ {\rm dist} (\rho) | \rho \in {\rm Corr} (M,N) \}  $ where ${\rm Corr} (M,N) := \{ A \subset M \times N | {\rm pr}_1 (A) = M, {\rm pr}_2 (A) = N\}$ is the set of correspondences (right-total and left-total relations) between $M$ and $N$ and ${\rm dist} (\rho) := \sup \{ |\tau_g (m_1, m_2) - \tau_h (n_1, n_2) | : (m_1, n_1), (m_2 , n_2) \in \rho \}$  is called the distortion of $\rho$, where for a globally hyperbolic $(Y,k)$, the map $\tau_k: Y \times Y \rightarrow [0; \infty)$ defined by $\tau_k(x,y) := \sup \{ \ell_k (c) | c: x \leadsto y {\rm \ causal} \} $ for $x \leq y$ (with $\ell_k$ the Lorentzian curve length) extended by $0$ on $Y \times Y \setminus J^+$ is called the {\bf Lorentzian distance function}. It has been shown in \cite{oM:LGH} that $d^-$ is a metric on the set of compact Cauchy slabs.

\newpage

Moreover, in the finite setting, we cannot assign, as before, to a Cauchy slab $X$ a subset of orders on $\N_K$ induced by some analogon to $V(X)$ any more (the condition on $V (X)$, being asymptotic, cannot be formulated any more on $\N_K$), instead, we apply a probabilistic construction: On the set $X^K$ of finite sequences $a: \N_K^* := \N \cap [1;K] \rightarrow X$ we define a map into the set of measures on order relations on $\N_K^*$ as follows: Let $K \in \N$, then $Q_K:= P(\N_K^* \times \N_K^*)$ is finite. Then

\[ M_K:=  \{f: Q_K \rightarrow [0; 1] \  \big|  \ \sum_{q \in Q_K } f(q) = 1 , f \ci (\s \times \sigma ) = f \ \forall \s: \N_K^* \hookrightarrow \N_K^* \}, \] 

i.e. $M_K$ is the set of probability measures on $Q_K$ invariant under permutations. For the product measure $\mu_K$ on $X^K$, define $C_K: CS \rightarrow M_K$,

\[ C_K (X) (q) := \mu_K (A_q (X)) , \qquad A_q (X):= \bigcup \{ A \subset X^K | a^*(\leq) = q \forall a \in A \}.\]

On the right-hand side of the map we consider {\em isomorphism classes} of elements of $P( \N_K^* \times \N_K^*) $: Let $S_K (U)$ be the space of $U$-preserving bijections\footnote{e.g. $S_K (\leq)= \{ \1_{\N_K^*}\}$ and $S_K(\1) = S_K$.} of $\N_K^*$, then the measure $PC_K(M) (U) := \sum_{I \in S_K (U)} \mu_K(I_* U ) = \# S_K (U)  $ is the probability for an order-preserving embedding for elements of the isomorphism class, we call it {\bf probability of the class}. $C_N (X,g)$ above is the pushforward $F_* (\vol_g^N)$ of $\vol_g^N$ by the pullback operator $F: X^{N+1} \ni a \mapsto a^* (\leq) \in {\rm Or} (\N_N)$.

\bigskip

The real number ${\rm tdiam} (Y) := \diam_{\R} (\tau_g(Y \times Y))$ is called {\bf timelike diameter} of $(Y,g)$, and it is not hard to see that {\rm tdiam} is increasing under inclusion or injective  maps preserving the Lorentzian distance, and that for $Y,Z \in CS$ we have $d^-(Y,Z) \geq | {\rm tdiam} (Y) - {\rm tdiam} (Z) |$. On $M_k \subset \R^{\N_k}$ we consider the total variation distance $ TV : (\mu, \nu ) \mapsto \sup_{A \subset M_k}  | \mu(A) - \nu (A)|$.

\begin{Theorem}[The finite Hauptvermutung is wrong for $d^-$]
	\label{HauptWrong}
Let $K \in \N$, $\e \in (0;1/2)$, $D>0$. Let $(X,g)$ any Cauchy slab with $\vol_g (X) =1$. Then there is $Y \in CS$ with ${\rm vol} (Y) =1$, $\partial^\pm X = \partial^\pm Y$, $d^-(X,Y) >D$ and $TV(C_K (X)  , C_K(Y) )< \e$.
\end{Theorem}

\V{\bf Proof.} Let $X $ be any Cauchy slab. There is $x \in X$ with $\vol_g (J^+(x)) < v/2$ with $(1-v)^K > 1- \e/2$. Let $c$ be a maximizer from $x$ to the future boundary of $X$. Let $u$ be a conformal factor with $u=1$ outside of a sufficiently thin neighborhood of $c$ s.t. $\vol_{u \cdot g} (J^+(x)) = \vol_g (J^+(x))$ and $\ell_{ug} (c)>{\rm tdiam} (X,g) +D $. Then, for $Y= (X, ug)$, we have $C_K(X,ug)(q) > \e$, and $d^- (X,Y) \geq {\rm tdiam} (X, ug) - {\rm tdiam} (X, g) > {\rm tdiam} (J^+(x), ug) - {\rm tdiam} (X, g) \geq  \ell(c, ug) - {\rm tdiam } (X,g)> D$. For the last assertion, first note that for the measures $\mu:= (\vol_g)^{\x N}$ and $\nu = (\vol_{ug})^{\x N}$ on $X^N$ we get $ |\mu (A) - \nu (A)| \leq |\mu (A \cap J^+(x)) + \mu (A \setminus J^=(x)) - (\nu (A \cap J^+(x) ) + \nu (A \setminus J^+(x))  | \leq |\mu (A \setminus J^+(x)) - \nu (A \setminus J^+(x) )| + | \mu (A \cap J^+(x) ) - \nu (A \cap J^+(x))|  $, and whereas the first term of that sum vanishes, the second is $\leq \e $. Then the claim follows from the classical fact that the $TV$ is decreasing (contractive) under pushforward of measures, and $C_K(X) = F_* (\vol_g)^N$. \hfill \qed

\bigskip

In some articles (e.g. not in \cite{Surya} but in \cite{DS}) the Hauptvermutung appears with additional hypotheses: In \cite{DS} it is required that 'The characteristic distance over which the continuum geometry $(M, g)$ varies appreciably is everywhere much greater than the Planck length/time'. However, it seems difficult to provide a rigorous Lorentzian-geometric reformulation of this requirement.

A requirement in \cite{DS} is called 'Planck-scale uniform' and defined by 'The number of causal set elements embedded in any sufficiently large, physically nice region of M is approximately equal to the spacetime volume of the region in fundamental, Planckian scale, volume units.'. It could be translated as the restriction to Planck-scale uniform finite sequences where, for $s>0$, a finite sequence $a: \N_K^* \rightarrow X$ into a Cauchy slab $X$ is called {\bf $s$-uniform} iff

\[ \forall U \subset X  {\rm \ closed \ causally \ convex \ }    \# (a^{-1} (\Int (U))) \in K \cdot (\vol (\Int(U)) - s  ; \vol(\Int(U)) + s ) , \]

and {\bf Planck-scale uniform} iff $n := \dim (X)$ and $s^{1/n} = h$ is the Planck length.

\begin{Theorem}
 For each sequence $a: \N\rightarrow X $ in a spatially compact Cauchy slab satisfying the volume law, there is $K \in \N$ such that $a|_{\N_K}$ is Planck-scale uniform.
 \end{Theorem}
 
\V{\bf Proof.} The set of causally convex subsets of $X$ is, by the limit curve theorem, a closed subset of the compact set of compact subsets of $X$, which is compact itself in the Hausdorff topology by Blaschke, and $s$-uniformity is an open condition. \hfill \qed
 
\bigskip

If $a$ is $s$-uniform, for every cone $J^\pm (x)$ and its complement, there is an inner and an outer approximation by some $J^\pm (a(B))$ where $B \subset \N_K$ with error at most $2s$ (just take $B:= \{ m \in \N | a(m) \leq x \}$ and observe that $C:= \{ m \in \N | a(m) \in J^-(a(B))\} =B$). 

\bigskip
 
This a priori requirement on $a$ could, rather unsactisfactorily, mean that fixing $K$ implies disregarding possibly physically relevant spacetimes, thus a crucial question in this context is: 

\medskip

{\em Given $D>0$, is there $K \in \N$ s.t. each $n$-dimensional Ricci-flat s.c. Cauchy slab $X$ with ${\rm tdiam} (X), \diam (\partial^- X), \diam (\partial^+ X) \leq D$ has a Planck-scale uniform finite sequence $a: \N_K^* \rightarrow X$?} 

\medskip

By scaling, it is quite easy to see that without diameter restrictions the assertion would be wrong.

The construction in Theorem \ref{HauptWrong} also shows that, under the additional requirement of Planck-scale uniformness alone, the Hauptvermutung for $d^-$ would still be wrong. 

However, we could pursue two mathematically sound loopholes to circumvent the result of Th. \ref{HauptWrong}:

\begin{enumerate}
	\item {\em Impose energy conditions on the class of admissible spacetimes}. We could hope that the modification above is less easily to perform if we wanted to satisfy energy conditions, too. This would still be a physically rather undesired loophole as the goal of CST is to model path-integral quantization including metrics and fields not critical for the Einstein-field Lagrangian.
	\item {\em Apply another notion of distance between Lorentzian length spaces}, e.g. one for measured Gromov-Hausdorff convergence (requiring that for each $\e>0$, a subset of measure $1- \e$ Gromov-Hausdorff converges), or the distance $d^\times $ on ordered measure spaces described below.   
\end{enumerate}

\newpage

Let us first exclude an realization of the first idea. Let $S(r)$ be the one-sphere of radius $r$ and let $C(T) = ((0, T) \times S(T^{-1/n})^n)$ be the flat normalized Lorentzian cylinder of timelike diameter $T$. A small calculation reveals: 

\begin{Theorem}[The finite vacuum Hauptvermutung is wrong for $d^-$]
For any $\e \in (0;1)$, $K \in \N$ there is $T>0$ such that for all $S>T$ the probability $E:= PC_{K} (C(S)) (R)$ of the isomorphism class $R$ of totally ordered relations on $\N_K$ satisfies  $E > \e  $. Furthermore, all $C_S$ are flat, in particular Einstein-vacuum solutions. However, $d^- (C_S, C_U ) \geq |{\rm tdiam} (C_S) - {\rm tdiam} (C_U)| = |S-U|$. 	
\end{Theorem}

\V{\bf Proof.} We calculate $\vol (X \setminus J(p)) \leq 4 \pi T^{-1/n}$ for all $p \in X$. As there are $K^2 $ different $2$-tuples of points, we get $E \geq 1 - K^2 \cdot 4 \pi T^{-1/n}$. \hfill \qed

\bigskip

Still, energy conditions plus diameter bounds could imply a version of the Hauptvermutung.

\bigskip

As to the second loophole, let ${\rm POM}_{\leq S}^I$ be the set of isomorphism classes of ordered measure spaces of cardinality $\leq S$. Furthermore, let ${\rm LBM} (\R) $ be the set of locally bounded measurable real functions on $\R $ and let $f \in {\rm LBM} (\R)$. To a point in an object $X$ of $CS$ we can assign the element $f \circ \tau_x$ of the space ${\rm AE}(X)$ of Borel-almost everywhere defined real functions on $X$ (where $\tau_x := \tau(x, \cdot)$ is the Lorentzian distance from $x$), then defining ($p \in [1; \infty]$ fixed) a map $\Phi_{f,p}$ takes the class of an object $(X, g)$ of $CS$ to the class of $ (X, d_{f,p})$ by

\vspace{-0.3cm}

\bea
\label{PhiDef}
\Phi_f : X \ni x \mapsto f \circ \tau_x, \qquad  d_{f,p}(x,y) := \Phi_{f,p} (\tau) := ((\Phi_f)^* d_{L^p}) (x,y)  = \vert f \ci \tau_x - f \ci \tau_y |_{L^p(X)}  \in [0; \infty ]
\eea

For $x \in X$, let $\tau_x^+$ resp. $\tau_x^-$ denote the positive resp. negative part of $\tau_x$. For $r \in [-1;1] $ we define

$$ F_r:= - (\frac{1}{2} - \frac{r}{2})\chi_{(- \infty; 0)} + (\frac{1}{2} + \frac{r}{2})\chi_{(0;  \infty)} : \R \mapsto \R  , \ D_r := d_{F_r,2}, $$

then $D_r$ interpolates between the past metric (taking into account only the past cones) $D_{-1} $ with 

$$D_{-1} (x,y) := \vert \vert  \chi_{(-\infty; 0)}\ci \tau_x -  \chi_{(-\infty; 0)} \ci  \tau_y \vert \vert_{L^2 (X)} = \sqrt{\mu (J^-(x) \triangle J^-(y))} $$

for $F_{-1} = (1- \theta_0)  $ and the future metric $D_1$ (taking into account only the future cones) for $F_1 = \theta_0 $, passing through $D_0 (x,y) = \frac{1}{2} \big || \sgn \tau_x - \sgn \tau_y \big||_{L^2}$. We define

$$\Phi^\times (X, \s) := (X, D_{-1/2}) \sqcup (X, D_0) \sqcup (X, D_{1/2})$$

and $d_{{\rm GH}}^\times := (\Phi^\times)^* d_ {{\rm GH}}$ where $d_{GH}$ is the usual Gromov-Hausdorff metric applied to the three-component metric spaces on the right-hand side. Let ${\rm POM^I_{\leq S, fv}}$ be the subset of ${\rm POM^I}_{\leq S}$ of all those classes s.t. the measure of future and past cones is finite. It was proven in \cite{oM:LGH} that $d^\times_{{\rm GH}} $ is an extended pseudometric on ${\rm POM^I_{\leq S}}$ and a metric on ${\rm POM^I_{\leq S, fv}}$. And the Malament theorem implies that $CS$ can be seen as a subset of $ {\rm POM}_{\leq L^2(\R), fv}$.

\begin{Theorem}[The Planck-scale uniform Hauptvermutung is true for $d^{\times} $]
	\label{HauptTrue}
	Let $s >0$. Then for each $\delta >0$ there is $K \in \N$ (actually, all $K > 8s/ \de $ work) such that for each two volume-normalized Cauchy slabs $X,Y$ such that there are $s$-uniform finite sequences $a: \N_K \rightarrow X$, $b: \N_K\rightarrow Y$ with $a^{-1} (\leq_X) = b^{-1} (\leq_Y)$ we get $d^\times (X,Y) < \delta$.   	
\end{Theorem}

\V {\bf Proof.} Let $X$ and $Y$ be as in the hypothesis of the theorem, then we define a correlation between $X$ and $Y$ as follows: For each $x \in X$ and each $y \in Y$ we define $x \sim y$ if and only if $a^{-1} (J^{\pm} (x)) = b^{-1} (J^{\pm} (y))$. Let $x_1, x_2 \in X$ be given and let $y_1, y_2 \in Y$ with $x_1 \sim y_1$ and $x_2 \sim y_2$. Let us first focus on $D_{-1/2}$, comparing the volume of the symmetric difference $\Delta (J^- (x_1), J^- (x_2))$ of the pasts of $x_1$ and $x_2$ to the same thing w.r.t. $y_1 $ and $y_2$. With $\Delta (J^-(x_1), J^-(x_2)) := (J^-(x_1 ) \setminus J^-(x_2)) \cup (J^-(x_2 ) \setminus J^-(x_1))$ (which is not causally convex if nonempty) and $P(u,v) := J^- (u) \cap J^- (v)$ (which is causally convex) we first state that whenever $x_1 \sim y_1$ and $x_2 \sim y_2$ then $a^{-1}( J^-(x_k)) = b^{-1} (J^-(y_k))$ and thus $|\vol (J^-(x_i)) - \vol (J^-(y_i))| < 2s/K $ and $ |\vol (P(x_1, x_2)) - \vol(P(y_1, y_2))| < 2s/K$. We calculate

\bean
|\vol (\Delta (J^- (x_1), J^- (x_2)))- \vol (\Delta (J^- (y_1, J^- (y_2)))) |\\
\leq |\vol (J^-(x_1) \setminus J^-(x_2)) - \vol (J^- (y_1) \setminus J^-(y_2))| + | \vol (J^-(x_2) \setminus J^-(x_1)) - \vol (J^- (y_2) \setminus J^-(y_1))  |\\
= |\vol (J^-(x_1)) - \vol (J^-(x_1) \cap J^-(x_2)) - \vol (J^- (y_1)) + \vol (J^- (y_1) \cap J^-(y_2))| \\
+ |\vol (J^-(x_2)) - \vol (J^-(x_2) \cap J^-(x_1)) - \vol (J^- (y_2)) + \vol (J^-(y_2) \cap J^- (y_1))|\\
\leq |\vol (J^- (x_1)) - \vol (J^- (y_1))| + | \vol (J^- (x_2) ) - \vol (J^- (y_2))| +  2|\vol (P(x_1,x_2)) - \vol (P(y_1, y_2))| \leq 8s/K,
\eean

so if $8s/K < \de$ the $D_{-1}$-distortion is $\leq \delta$. Similar estimates hold for $D_0$ and $D_{\pm 1/2}$. \hfill \qed

\bigskip

At a first sight, the difference between $d^-$ and $d^\times$ apparent in Theorems \ref{HauptWrong} (and the persistencee of its result even under the assumption of Planck-scale uniformness) and \ref{HauptTrue} may surprise a bit, considering that in \cite{oM:LGH} it has been shown that $d^\times$ and $d^-$ generate the same topology on the closure of the class $CS $ of Cauchy slabs w.r.t. either metric, but this does not mean that they generate the same uniformity on $CS$ --- in fact, Theorems \ref{HauptWrong} and \ref{HauptTrue} show that they do not.

More refined statements can be obtained if we define, for a Cauchy slab $X$, the space $LPC_{K,s} (X)$ of measures on the order relations on $\N_K$ induced by $s$-uniform finite sequences $a: \N_K \rightarrow X$.

\bigskip

\V $CS$ can be completed with morphisms (time-oriented isometric embeddings with causally convex image) to obtain a category, and the same is true for  ${\rm POM}$ (here we take ordered measure-preserving embeddings with causally convex image). On the range of the map $C$, the morphisms we can consider are the following: For $A,B \in {\rm or} (\N)$, a morphism $F$ from $A$ to $B$ is a choice, for each $a \in A$, of some $b \in B$ and of an increasing map $(\N, a) \to (\N, b)$. Then from a morphism $f: M \to N$ in the category $CS$ (or else ${POM}_{\leq S, fv}$) we can define a morphism in the category ${\rm or} (\N)$ by $a^* \leq_M \mapsto (f \ci a )^* \leq_N$, and it is easy to see that this is well-defined and defines a functor between the categories.

\bigskip

\bigskip

\newpage

The author wants to acknowledge helpful comments of Faye Dowker, Jona Röhrig, Sumati Surya, Johan Nedin and of two anonymous referees on previous versions of this paper.

\medskip

This research was funded in part by the Austrian Science Fund (FWF) [Grant DOI 10.55776/EFP6].

\begin{small}

\end{small}

\end{document}